\documentclass[a4paper,11pt,centertags]{amsart}
\usepackage{xspace}
\usepackage{a4wide}

\newcommand\R{\ensuremath{\mathbb{R}}\xspace}  
\newcommand\N{\ensuremath{\mathbb{N}}\xspace} 
\newcommand{\aop}{{\mathbf a}}  
\newcommand\dleb{\operatorname{d}}
\newcommand\dx{\ensuremath{\dleb\!x}}

\newcommand\diw{\operatorname{div}}
\newcommand\thec{\ensuremath{\hat{\theta}}\xspace}
\newcommand\ucha{\ensuremath{\hat{u}}\xspace}
\newcommand\uep{\ensuremath{u^\varepsilon}}
\newcommand\thep{\ensuremath{\theta^\varepsilon}}
\newcommand\fep{\ensuremath{f^\varepsilon}}
\newcommand\ftid{\ensuremath{\widetilde{f}}}
\newcommand\A{{\mathbf A}}

\newcommand\indi{\nbOne}
\newcommand\nbOne{{\mathchoice {\mathrm {1\mskip-4.1mu l}} {\mathrm{ 1\mskip-4.1mu
        l}} {\mathrm {1\mskip-4.6mu l}} {\mathrm {1 \mskip-5.2mu l}}}}

\makeatletter
\def\th@definition{%
  \normalfont 
}
\def\th@plain{%
  \slshape 
}
\def\th@remark{%
  \normalfont 
  \thm@preskip\topsep
  \divide\thm@preskip\tw@
  \thm@postskip\thm@preskip
}
\makeatother
\theoremstyle{remark}
\newtheorem{rmk}{Remark}[section]

\theoremstyle{definition}
\newtheorem{defi}[rmk]{Definition}
\theoremstyle{plain}
\newtheorem{prop}[rmk]{Proposition}
\newtheorem{lemma}[rmk]{Lemma}
\newtheorem{theorem}[rmk]{Theorem}
\newtheorem{coro}[rmk]{Corollary}

\author[{O. Guib\smash{\'e}}]
{Olivier Guib\'e}

\title{Existence and uniqueness results 
for a nonlinear stationary system}

\keywords{nonlinear coupled system,
renormalized solutions, $L^1$ data}

\begin{document}

\maketitle
\ \\ \vglue-1.2cm
\hglue0.02\linewidth\begin{minipage}{0.9\linewidth}
\begin{center}
{Laboratoire de~Math\'ematiques~Rapha\"el~Salem}\\ 
UMR CNRS 6085, Site Colbert \\
Universit\'e de Rouen\\
F-76821 Mont Saint Aignan cedex \\
 E-mail : \texttt{olivier.guibe@univ-rouen.fr} 
\end{center}
\end{minipage}

\vglue 0.5cm

\normalsize
\thispagestyle{empty}
\medskip

\begin{abstract}
We prove a few existence results of a solution for a static system
with a coupling of thermoviscoelastic type. As this system involves $L^1$ 
coupling terms we use the techniques of renormalized solutions for
elliptic equations with $L^1$ data. We also prove partial uniqueness results.
\end{abstract}

\medskip

\section{Introduction}
In the present paper we consider the following nonlinear coupled system:
\begin{align}
  \ & \lambda u - \diw\big( \A(x)Du-f(\theta)\big) =
 g\quad  & & \text{ in $\Omega$,} \label{eq1} \\
\ &  \mu \theta -\diw\big(\aop(x,D\theta)\big)= \big(\A(x)Du-
f(\theta)\big)\cdot
Du  & & \text{ in $\Omega$,} \label{eq2} \\
& u=0\qquad\theta=0 & &\text{ on $\partial\Omega$,} \label{eq3}
\end{align} 
where $\Omega$ is an open and bounded subset of $\R^N$ ($N\geq 2$),
$\lambda,\mu>0$, $f\,:\,\R\mapsto \R^N $ is a continuous function,
$g\in L^2(\Omega)$, $\A(x)$ is a coercive matrix with
$L^\infty$--coefficients and $v\mapsto-\diw\big(\aop(x,Dv)\big)$ is a
monotone operator defined from $H^1_0(\Omega)$ into $H^{-1}(\Omega)$.
\par
\smallskip
Problem (\ref{eq1})--(\ref{eq3}) 
is a static version (or time discretized-version) of a time dependent
class of systems in thermoviscoelasticity. Under stronger assumptions
than in the present paper, existence of a solution for these evolution
systems is established in \cite{BlGu99} (see also
\cite{BlGu97}). Moreover, for (\ref{eq1})--(\ref{eq3}) 
uniqueness results were also proved.
\par
The main difficulties in dealing with existence of solution of system
(\ref{eq1})--(\ref{eq3}) are due to equation (\ref{eq2}) and the coupling. Indeed if $u$ 
is a variational solution of (\ref{eq1}) (i.e. $u\in H^1_0(\Omega)$) then
the right--hand side of (\ref{eq2}) belongs to  $L^1(\Omega)$.
It follows from  L. Boccardo and T. Gallou\"et 
\cite{BG89} (see also \cite{BBGGPV95} and \cite{M94}) that $\theta$ is
expected in $L^q(\Omega)$ for $q<N/(N-2)$ if $N\geq 3$ and $q<\infty$
if $N=2$.
With the aim of solving (\ref{eq1}) with $f(\theta)\in L^2(\Omega)$ we 
are then led to assume that $f$ satisfies the growth assumption
\begin{equation*}
  \forall r\in\R\quad \big|f(r)\big|\leq a +M|r|^\alpha,
\end{equation*}
with $a>0$, $M>0$ and $\alpha<N/\big(2(N-2)\big)$ if $N\leq 3$ and
$\alpha<\infty$ if $N=2$. 
Under this hypothesis on $f$, the coupling between Equations
(\ref{eq1}) and (\ref{eq2}) together with the L. Boccardo and
T. Gallou\"et estimates techniques (see \cite{BG89} and 
Remark \ref{ogrmk55} of the present paper) lead to the following a priori
estimate on $\theta$,
\begin{equation*}
  \forall q\in[1,\frac{N}{N-2}[\quad\quad \|\theta\|_{L^q(\Omega)}
  \leq C\big( 1+ \|f(\theta)\|^2_{L^2(\Omega)}\big) \leq C'\Big(1+\int_\Omega |\theta|^{2\alpha} \dx \Big),
\end{equation*}
where $C$ and $C'$ are real positive constant independ of $u$ and
$\theta$. This implies that
if $2\alpha\geq 1$ the estimate above
is not sufficient in general settings to obtain
the existence of a solution of (\ref{eq1})--(\ref{eq3}) using a
fixed--point or approximation method.

\par

As the right-hand side of (\ref{eq2}) belongs to $L^1$,
 we use in the present paper the convenient framework of renormalized
solutions that insures uniqueness and stability results for equations 
with $L^1$
data. Renormalized solutions have been introduced by R.J.~DiPerna and
P.-L.~Lions in \cite{DPL89a} and \cite{DPL89b} for first order
equations and have been adapted for elliptic equations
in \cite{BDGM93}, \cite{LM}, \cite{M93} and for elliptic equations
with general measure data in \cite{DMOP99} (see also \cite{DMOP97}).
 Other frameworks as entropy solutions \cite{BBGGPV95} or 
SOLA \cite{Dall96} may be used for equation (\ref{eq2}) with $L^1$
data.
\par
Another interesting question related to problem
(\ref{eq1})--(\ref{eq3}) deals with the uniqueness of a solution, that 
is an open problem in general settings due to lack of regularity of
$\theta$ and the right--hand side of (\ref{eq2}). We investigate in
the present paper  uniqueness of a small solution
$(u,\theta)$ such that $\theta\geq 0$ almost everywhere in $\Omega$ and under
additional assumptions on the data for $N=2$ and $N=3$.
\par
Elliptic systems involving $L^1$ coupling terms are also studied
in \cite{GH}, \cite{CFC97} and \cite{Lew97} and use a convenient
formulation for equation with the $L^1$ term.
\par
\smallskip
The plan of this paper is as follows. In Section 2 we recall the
definition of a renormalized solution and we define a weak-renormalized 
solution for system (\ref{eq1})--(\ref{eq3}). In Section 3 we give a
few useful properties on  renormalized solutions.
Section 4 and Section 5 are devoted to existence results for two
restricted case: the first case deals with small data, the second case 
contains existence results under more restrictive conditions on $f$ but 
for general data.
Section 6 contains a partial uniqueness result of a small solution
$(u,\theta)$ such that $\theta\geq 0$ almost everywhere in $\Omega$ and under
additional assumptions on the data.
\par\smallskip

\section{Assumptions and definitions}
Let $\Omega$ be an open and bounded subset of $\R^N$ ($N\geq 2$). The
following assumptions are made on the data:
\renewcommand{\labelenumi}{(A\theenumi)}
\makeatletter\renewcommand{\p@enumi}{A}\makeatother
\begin{enumerate}
\item\label{as1} $\A(x)$ is a coercive matrix field with coefficients lying in
  $L^\infty(\Omega)$ i.e.
  $\A(x)=\bigl(a_{i,j}(x)\bigr)_{1\leq i,j\leq N}$ with \\
  \ $\bullet$  $a_{i,j}(x)\in L^\infty(\Omega)$\\
  \ $\bullet$ $\exists \gamma>0 $ such that $\forall \xi\in\R^N$
  $\A(x)\xi\cdot \xi\geq \gamma\|\xi\|^2_{\R^N}$\quad for almost
  every $x\in\Omega$;
\item\label{as2a} the function $\aop\ :\ \Omega\times\R^N\longmapsto\R^N$ is a
  Caratheodory function (continuous in $\xi$ for almost every $x\in\Omega$ and
  measurable in $x$ for every $\xi\in\R^N$) and there exists $\delta>0$ such
  that
  $$\forall \xi\in\R^N\quad \aop(x,\xi)\cdot \xi \geq \delta
  \|\xi\|_{\R^N}^2 \quad \text{for almost every $x\in\Omega$;}$$
\item\label{as2b} for every $\xi$ and $\xi^\prime$ in $\R^N$, 
  and almost everywhere in $\Omega$
  $$(\aop(x,\xi)-\aop(x,\xi^\prime))\cdot(\xi-\xi^\prime)\geq 0 ;$$
\item\label{as2c} there exists $\beta>0$ such that
  $$
  |\aop(x,\xi)|\leq\beta(b(x)+|\xi|)$$
  holds for every $\xi\in\R^N$
  and for almost every in $x\in\Omega$ with $b\in L^2(\Omega)$;
\item\label{as3}  $\lambda>0$, $\mu> 0$;
\item\label{as4} $f$ is a continuous function defined on $\R$ with values in
  $\R^N$ ;
\item\label{as5} $g$ is an element of $L^2(\Omega)$.
\end{enumerate}

Throughout this paper and for any non negative real number $K$ we
denote by $T_K(r)$ the truncation function at height $\pm K$, i.e.
$T_K(r)=\min\bigl(K,\max(r,-K)\bigr)$.
 For a measurable set $E$ of $\Omega$, we denote by
$\indi_{E}$ the characteristic function of $E$.

Following \cite{LM} (and \cite{M93})
 we recall the definition of a renormalized solution for 
nonlinear equations of type (\ref{eq2}) with $L^1$ right--hand side.
\begin{defi}\label{defsr}
  Let $F$ be an element of $L^1(\Omega)$.  A measurable function
  $\theta$ defined on $\Omega$ is called a renormalized solution of
  the problem
  $$P(F)\left\{\begin{aligned}
      \ & \mu\theta-\diw (\aop(x,D\theta))=F\ \ \text{ in }\Omega, \\
      \ & \theta=0 \ \ \text { on }\partial\Omega,
\end{aligned}\right.$$ 
if
\begin{align}
  \ & \theta\in L^1(\Omega),\  \forall K>0 \ \ 
T_K(\theta)\in H^1_0(\Omega); \label{defsr1}\\
  \text{for }&\text{every function $h\in W^{1,\infty}(\R)$ such that
    $h$ has a compact support, } \notag \\
  \ &\mu\theta
  h(\theta)-\diw\bigl(h(\theta)\aop(x,D\theta)\bigr)
+h^\prime(\theta)\aop(x,D\theta)\cdot
  D\theta=Fh(\theta)
  \text{ in }{\mathcal D}^\prime(\Omega); \label{defsr2}\\
  \ &\lim_{n\rightarrow\infty}\frac{1}{n}\int_{\{n<|\theta|<2n\}}
  |D\theta|^2 \dx =0. \label{defsr3}
\end{align}
\end{defi}

  Under assumptions (\ref{as2a})--(\ref{as2c}) and with $\mu>0$,
 using the techniques developped in \cite{LM} (see also \cite{DMOP99}
 and \cite{M93}), 
there exists a unique renormalized solution of $P(F)$ for any $F$
in $L^1(\Omega)$.
\par
We now use renormalized solutions 
to define a so called weak-renormalized 
solution of Problem (\ref{eq1})--(\ref{eq3}).

\begin{defi}\label{defswr}
A couple of functions $(u,\theta)$ defined on $\Omega$ is called a
weak-renormalized solution of (\ref{eq1})--(\ref{eq3}) if $u$ and
$\theta$ satisfy
\begin{gather}
  u\in H^1_0(\Omega), \label{defswr1} \\
f(\theta)\in \big(L^2(\theta)\big)^N ,\label{defswr2} \\
\lambda u -\diw\big( \A(x)Du-f(\theta)\big)=g \quad\text{in }{\mathcal
  D}'(\Omega), \label{defswr3} \\
\text{$\theta$ is a renormalized solution of (\ref{eq2})--(\ref{eq3}).}
\label{defswr4}
\end{gather}
\end{defi}
\par\medskip

Under regularities (\ref{defswr1})--(\ref{defswr2}), the right--hand 
side of (\ref{eq2}), $(\A(x)Du-f(\theta))\cdot Du$, belongs to
$L^1(\theta)$. So we are in the framework of renormalized solution for
equation (\ref{eq2}).

\section{Useful properties of renormalized solutions}

We recall the following propositions on  renormalized solutions of
elliptic equations for $L^1$ data, that can be shown using the
techniques developped in \cite{DMOP99}, \cite{LM} and
\cite{M93}.

\begin{prop}[Existence and uniqueness of the renormalized solution]
\label{prop1}
Assume that  (\ref{as2a})--(\ref{as2c}) hold true and $\mu >0$.
Then for any $F$
belonging to $L^1(\Omega)$, there exists a unique renormalized
solution of Problem $P(F)$. Moreover 
for any function $w\in L^\infty(\Omega)\cap H^1_0(\Omega)$, if 
there exists $K>0$ such that $Dw=0$ almost everywhere in $\{x\,:\,|\theta(x)|\geq
K\}$ then
\begin{equation}\label{p1p1}
\mu \int_\Omega \theta w\dx +\int_\Omega \aop(x,D\theta)\cdot Dw \dx
=\int_\Omega Fw\dx.
\end{equation}
\end{prop}

\begin{rmk}
  Equality (\ref{p1p1}) which is proved in \cite{DMOP99} in the
  context of general measure data, is formally obtained through using the
  test function  $w$ in the equation of $P(F)$.
\end{rmk}

\begin{prop}\label{prop2}
Assume that (\ref{as2a})--(\ref{as2c}) 
hold true and $\mu>0$. Let $F_1$, $F_2$ be
two elements of $L^1(\Omega)$, and denote by $\theta_i$ the unique
renormalized solution of $P(F_i)$ $(i=1,2)$.
\par
Then for any $K>0$
\begin{multline}
  \label{p2p1}
  \mu \int_\Omega (\theta_1-\theta_2) T_K(\theta_1-\theta_2) \dx
+ \int_{\{|\theta_1-\theta_2|<K\}} \big(\aop(x,D\theta_1)-\aop(x,D\theta_2)\big)\cdot
(D\theta_1-D\theta_2) \dx \\
\leq \int_\Omega (F_1-F_2) T_K(\theta_1-\theta_2)\dx.
\end{multline}
\end{prop}

\begin{rmk}
  \label{rmk1} Inequality (\ref{p2p1}) is obtained by plugging the
  admissible test function \linebreak
$h_n(\theta_1) h_n(\theta_2)
  T_K(\theta_1-\theta_2)$ in the difference of the equations $P(F_1)$
  and $P(F_2)$ (that is licit in view of  Proposition
  \ref{prop1}) where $h_n$ is a sequence of functions in
  $W^{1,\infty}(\R)$ such that $h_n(r)\rightarrow 1$ as $n$ tends to
  $\infty$ and with compact support.
\par
Due to Proposition \ref{prop2} we deduce that 
$$\big\|\theta_1-\theta_2\big\|_{L^1(\Omega)} \leq \frac{1}{\mu}
\big\|F_1-F_2\big\|_{L^1(\Omega)}
$$
and the continuity of the renormalized solution of $P(F)$ 
with respect to the datum $F$.
\end{rmk}

We recall the following lemma  that can be proved by means of the
estimates techniques
of L.~Boccardo and T. Gallou\"et \cite{BG89} (see
also \cite{BBGGPV95}).

\begin{lemma}\label{lem1}
Let $\theta$ be a measurable function defined on $\Omega$, that is
finite almost everywhere in $\Omega$,  and $M>0$ 
such that 
\begin{gather*}
  \forall K>0,\quad T_K(\theta)\in H^1_0(\Omega)\quad \text{ and } \\
\int_\Omega \big|DT_K(\theta)\big|^2\dx \leq K M.
\end{gather*}
Then $\theta\in W^{1,p}_0(\Omega)$  for any  $1\leq p <N/(N-1)$  and
there exists a constant $C$ (depending on $\Omega$ and $p$) such that
$$\big\|\theta\big\|_{W^{1,p}_0(\Omega)} \leq C M.
$$  
\end{lemma}

Gathering Proposition \ref{prop1} and Lemma \ref{lem1} we deduce the
following corollary.

\begin{coro}\label{coro1}
  Assume that (\ref{as2a})--(\ref{as2c}) hold true and
  $\mu>0$. Let $F$ be an element of $L^1(\Omega)$ and $\theta$ the 
  renormalized solution of $P(F)$. Then for any $1\leq p < N/(N-1)$,
  $\theta\in W^{1,p}_0(\Omega)$ and
  \begin{equation}
    \label{cor1}
    \|\theta\|_{W^{1,p}_0(\Omega)} \leq C \|F\|_{L^1(\Omega)}
  \end{equation}
where $C$ is a constant only depending upon $\Omega$, $p$ and $\aop$.
\end{coro}


\section{Existence of small solutions of (\ref{eq1})--(\ref{eq3}) for
  small data}
In this section we assume that the continuous function $f$ satisfies
the following growth assumption
\begin{equation}
  \label{sm1}
  \exists a\geq 0,\exists M>0,\quad \forall r\in\R\quad
\big|f(r)\big|\leq a+M|r|^\alpha 
\end{equation}
with $1/2<\alpha < N/(2(N-2))$ if $N\geq 3$ and $1/2<\alpha<+\infty$
if $N=2$.
\par
Under this additional assumption, Theorem \ref{th1} insures the existence of at
least a solution of Problem (\ref{eq1})--(\ref{eq3}) for small data
enough.
Notice that on the one hand the upper bound of $\alpha$ in (\ref{sm1})
is motivated in Introduction, on the other hand the lower bound
permits us to exploit the small character on the data.

\begin{theorem}\label{th1}
  Assume that (\ref{as1})--(\ref{as5}) and (\ref{sm1}) hold
  true. There exists a real positive number $\eta$ such that if
  $a+\|g\|_{L^2(\Omega)} <\eta$, then there exists at least a
  weak-renormalized solution of (\ref{eq1})--(\ref{eq3}) such that
  \begin{equation*}
    \label{sm2}
    \big\|u\big\|_{H^1_0(\Omega)} +
    \big\|\theta\big\|_{L^{2\alpha}(\Omega)} \leq \omega(\eta)
  \end{equation*}
where $\omega(\eta)$ tends to zero as $\eta$ tends to zero.
\end{theorem}

\begin{proof}
 The proof is divided into 2 steps. Step 1 is devoted to the
 construction of a fixed--point operator. In Step 2 we give a
 sufficient condition on the data in order to apply  the Schauder fixed--point Theorem.
\par
\smallskip
{\noindent \slshape Step 1.} Since the function $f$ is continuous and
verifies growth assumption (\ref{sm1}), under  assumptions (\ref{as1}),
(\ref{as3}) and (\ref{as5}) the mapping 
\begin{align}
L^{2\alpha}(\Omega)  & \longmapsto   H^1_0(\Omega) \notag \\
\thec & \longmapsto \ucha,\quad \text{where $\ucha$ is the unique
  solution of the problem } 
\notag \\
\null  &\qquad\quad
\left\{
\begin{gathered} 
\lambda \ucha - \diw\big( \A(x)D\ucha-f(\thec)\big) =
 g\quad \text{ in $\Omega$}  \\
\ucha=0 \quad \text{ on $\partial\Omega$,}
\end{gathered}
\right.
 \label{psm2}
\end{align}
is continuous and the coercivity of $\A$ implies that
\begin{equation}
  \label{psm4}
  \int_\Omega\big|D\ucha\big|^2\dx \leq C \Big(
  \big\|g\big\|_{L^2(\Omega)}^2 + 
\big\|f(\thec)\big\|_{(L^2(\Omega)^N)}^2 \Big)
\end{equation}
where $C$ is a generic constant independent of $\thec$.
\par
Let $\thec$ be an element of $L^{2\alpha}(\Omega)$ and $\ucha$ the unique
element of $H^1_0(\Omega)$ solution of
(\ref{psm2}). Due to growth assumption (\ref{sm1}) on
$f$ and the regularity of $\ucha$, the field $(\A(x)D\ucha-f(\thec))\cdot
D\ucha$ belongs to $L^1(\Omega)$, and by Proposition \ref{prop1}, 
let $\theta$ be the unique
renormalized solution of the problem:
\begin{align}
\ &  \mu \theta -\diw\big(\aop(x,D\theta)\big)= \big(\A(x)D\ucha-f(\thec)\big)\cdot
D\ucha  & & \text{ in $\Omega$} \label{sm5} \\
& \theta=0 & &\text{ on $\partial\Omega$.} \label{sm6}
\end{align}
We denote by $\Gamma$ the mapping defined by $\theta=\Gamma(\thec)$.
\par
Since $1<2\alpha < N/(N-2)$ (and $1< 2\alpha <+\infty$ if $N=2$), let
$q$ be a positive real number such  that $2\alpha< q^* < N/(N-2)$ (and $2\alpha <q^*
<+\infty$ if $N=2$), where $q^*$ denotes the Sobolev conjugate
exponent ($1/q^*=1/q -1/N$).\par
Using the properties of the renormalized solutions (see Remark
\ref{rmk1} and Corollary \ref{coro1}), the interpolation of $L^{2\alpha}(\Omega)$
between $L^1(\Omega)$ and $L^{q^*}(\Omega)$ and 
the Rellich Kondrachov Theorem we
deduce that $\Gamma$ is defined continuous and compact 
from $L^{2\alpha}(\Omega)$ into itself. Moreover inequality
(\ref{psm4}) and Corollary \ref{coro1} imply that $\forall \thec\in
L^{2\alpha}(\Omega)$, if $\theta=\Gamma(\thec)$ then
\begin{equation*}
  \big\|\theta\big\|_{W^{1,q}_0(\Omega)} \leq C \Big(
  \big\|g\big\|_{L^2(\Omega)}^2 +
  \big\|f(\thec)\big\|_{(L^2(\Omega))^N}^2\Big) 
\end{equation*}
and growth assumption (\ref{sm1}) on $f$ yields
\begin{equation}\label{psm7}
  \big\|\theta\big\|_{W^{1,q}_0(\Omega)} \leq C \Big(
  \big\|g\big\|_{L^2(\Omega)}^2 +a^2
+M^2 \big\|\thec\big\|_{L^{2\alpha}(\Omega)}^{2\alpha}\Big) 
\end{equation}
where $C$ is a constant independent of $\thec$.
\par
\medskip
{\noindent\slshape Step 2.} Applying the Schauder fixed--point Theorem to
the mapping $\Gamma$ reduces to show that there exists for instance a ball $B$ of
$L^{2\alpha}(\Omega)$ such that $\Gamma(B)\subset B$.
\par
Let $\thec$ be an element of $L^{2\alpha}(\Omega)$ and
$\theta=\Gamma(\thec)$. Since $1<2\alpha< q^*<N/(N-2)$ (and
$1<2\alpha<q^*<+\infty$ if $N=2$), the Sobolev embedding Theorem and
(\ref{psm7}) lead to
\begin{equation}
  \label{psm8}
  \big\|\theta\big\|_{L^{2\alpha}(\Omega)} \leq C\Big(
 \big\|g\big\|_{L^2(\Omega)}^2 +a^2
+M^2 \big\|\thec\big\|_{L^{2\alpha}(\Omega)}^{2\alpha}\Big)
\end{equation}
where $C$ is a constant independent of $\thec$, $g$, $a$ and $M$.
\par
As $2\alpha>1$, let $\eta>0$ and $R(\eta)>0$ such that
\begin{gather*}
  C\big(\eta+M^2 (R(\eta))^{2\alpha}\big)< R(\eta), \\ 
R(\eta) < 2C\eta. 
\end{gather*}
If
\begin{equation}
  \label{psm13}
  \big\|g\big\|_{L^2(\Omega)}^2 + a^2 < \eta
\end{equation}
then we have
\begin{equation*}
  \Gamma\Big(B_{L^{2\alpha}(\Omega)}\big(0,R(\eta)\big)\Big)
\subset B_{L^{2\alpha}(\Omega)}\big(0,R(\eta)\big).
  \end{equation*}
Therefore, we may apply the  Schauder fixed--point Theorem so that,
 there exists at least a solution $(u,\theta)$ of
(\ref{eq1})--(\ref{eq3}) in the sense of Definition \ref{defswr}.
\par
Moreover the choice of $R(\eta)$ and (\ref{psm4})
 imply that 
 \begin{equation*}
    \big\|u\big\|_{H^1_0(\Omega)} +
    \big\|\theta\big\|_{L^{2\alpha}(\Omega)} \leq \omega(\eta)
  \end{equation*}
where $\omega(\eta)$ tends to zero as $\eta$ tends to zero.
\par
The proof of Theorem \ref{th1} is complete.
\end{proof}

\section{Existence of a solution of (\ref{eq1})--(\ref{eq3}) for more
  general data}
In order to remove the small character on the data of the previous
section,
we suppose by now more restrictive hypotheses on the behavior of $f$,
which are on $\R^+$ 
\begin{equation}
  \label{g1}
  \left\{ 
    \begin{aligned}[c]
      \ & \lim_{r\rightarrow+\infty} \frac{|f(r)|}{r^{(N+2)/(2N)}} = 0
      \quad\text{ if $N\geq 3$}, \\
 &\forall r\in\R^+ \quad
 |f(r)|\leq a+ M|r| \quad \text{ if $N=2$ with $a\geq 0$ and $M\geq0$,}
    \end{aligned}\right.
\end{equation}
that is more restrictive than (\ref{sm1}) since
$(N+2)/(2N)<N/(2(N-2))$),
and  on $\R^-$ a behavior
of which the  model case is $f=0$ for $r<r_0\leq0$. (see hypotheses (\ref{g2}) and
(\ref{g3}) below).

\begin{theorem}\label{th2}
 Assume that assumptions (\ref{as1})--(\ref{as5}) and
  (\ref{g1}) hold true. Moreover assume that the continuous function
  $f$ is such that
  \begin{equation}
    \label{g2}
    \lim_{r\rightarrow-\infty} \frac{|f(r)|}{\sqrt{|r|}}=0.
  \end{equation}
Then there exists at least a weak-renormalized solution of
(\ref{eq1})--(\ref{eq3}).
\end{theorem}

In the case where the function $f$ has a zero  on $\R^-$, the structure
of equation (\ref{eq2}) allows us to remove the growth assumption on $f$ on $\R^-$
and give an additional property on $\theta$.

\begin{theorem}\label{th3}
  Assume that assumptions (\ref{as1})--(\ref{as5}) and
  (\ref{g1}) hold true. Moreover assume that the continuous function
  $f$ satisfies:
  \begin{equation}
    \label{g3}
    \exists r_0\in\R^-\quad\text{ such that}\quad f(r_0)=0.
  \end{equation}
Then there exists at least a weak-renormalized solution $(u,\theta)$ of
(\ref{eq1})--(\ref{eq3}) such that $\theta\geq r_0$ almost everywhere  in $\Omega$.
\end{theorem}

\begin{rmk}
The existence results of Theorems \ref{th2} and \ref{th3} were
  announced in \cite{Gu98} (see also  \cite{OGthese}) under more
  restrictive hypotheses on the function $f$. Let us notice that when $N=2$
  linear growth on $\R^+$ is allowed for $f$.
\end{rmk}

Before to prove Theorem \ref{th2}, we give a technical lemma.
\begin{lemma}\label{lem2}
 Assume that (\ref{g1}) holds true.
  Let $\theta$ be a measurable function defined on $\Omega$ such that
  \begin{gather}
    \theta\in L^1(\Omega), \label{l1} \\
\forall K>0\quad T_K(\theta)\in H^1_0(\Omega), \notag \\
\exists C_1>0\text{ such that }\forall K>0\quad\ 
\frac{1}{K} \int_\Omega \big|DT_K(\theta)\big|^2\dx < C_1 \Big(
\int_{\{0\leq \theta \leq K\}} f^2(\theta)\dx +1\Big).
\label{l2}
  \end{gather}
Then for any $1\leq q < N/(N-2)$ (and $1\leq q <+\infty$ if $N=2$),
there exists a constant $C'$, only depending upon $q$, $\Omega$,
$\|\theta\|_{L^1(\Omega)}$, $C_1$  and $f$ such that
\begin{gather}\label{len2}
\big\|\theta\big\|_{L^q(\Omega)} \leq C'.
\end{gather}
\end{lemma}

\begin{proof}[Sketch of the proof]
The proof relies on  estimate techniques  of L.~Boccardo and
T.~Gallou\"et \cite{BG89} (see also \cite{BBGGPV95}). If
$N=2$ we use the limit case of the Sobolev embedding Theorem
(see \cite{Adams}, \cite{GiTr} for instance) that allows us to reach linear
growth on $\R^+$ for the function $f$.
\par\smallskip
\noindent{\slshape Case $N\geq3$.}
 Let $n$ be an element of $\N$, that will
fixed in the sequel, and let $q$ be such that $1<q<N/(N-2)$. 
Hypothesis (\ref{g1}) gives
that
\begin{equation}
  \label{l1n1}
  \forall r\in\big[2^{n},+\infty\big[\quad
|f(r)|\leq \omega(n) r^{(N+2)/2N},
\end{equation}
where $\omega(n)$ tends to zero as $n$ tends to infinity.
\par
 As $\theta$ is finite almost everywhere  in $\Omega$, we  have
\begin{align*}
  \int_\Omega |\theta|^q\dx & \leq 2^{nq}|\Omega|
+\sum_{k=n}^{+\infty} \int_{\{2^k<|\theta|\leq 2^{k+1}\}} |\theta|^q\dx
\\
& \leq  2^{nq}|\Omega| + \sum_{k=n}^{+\infty}
\left(\frac{1}{2^k}\right)^{2^*-q} \int_\Omega
|T_{2^{k+1}}(\theta)|^{2^*}\dx,
\end{align*}
where $2^*$ denotes the Sobolev conjugate exponent
($1/{2^*}=1/2-1/N$). 
\par
The Sobolev embedding Theorem and (\ref{l2}) (with $K=2^{k+1}$) yield
\begin{gather*}
   \int_\Omega |\theta|^q\dx  \leq 2^{nq}|\Omega|
+ CC_1^{2^*/2}\sum_{k=n}^{+\infty} 
\left(\frac{1}{2^k}\right)^{2^*-q} \left( 2^{k+1} \int_{\{0\leq \theta\leq
    2^{k+1}\}} f^2(\theta)\dx + 2^{k+1}\right)^{2^*/2},
\end{gather*}
where $C$ is a constant depending on $\Omega$.
\par
Using (\ref{l1n1}) we obtain 
\begin{multline}
\label{l1n1a}
     \int_\Omega |\theta|^q\dx   \leq 2^{nq}|\Omega|
+  C(2C_1)^{2^*/2}\sum_{k=n}^{+\infty} 
\left(\frac{1}{2^k}\right)^{2^*/2-q} \Bigg(
  |\Omega|\max_{r\in [0,2^n]}|f(r)|^2+1 
\\
+\omega(n)\int_{\{2^n\leq \theta\leq 2^{k+1}\}}
|\theta|^{(N+2)/N}\dx\Bigg)^{2^*/2}.
\end{multline}

\par
On the one hand
H\"olders inequality gives, $\forall n\leq k <+\infty$,
\begin{equation}
  \label{l1n1b}
  \int_{\{2^n\leq \theta\leq 2^{k+1}\}}
|\theta|^{(N+2)/N}\dx \leq
\Big(\|\theta\|_{L^1(\Omega)}\Big)^{2/N}\,
\left(\int_{\{2^n\leq \theta\leq 2^{k+1}\}} |\theta|^{N/(N-2)}\dx
\right)^{2/2^*},
\end{equation}
on the other hand as $q<N/(N-2)=2^*/2$, the series $\sum_{k=n}^{+\infty}
\big(\frac{1}{2^k}\big)^{2^*/2-q}$ is convergent and we have
\begin{gather}\label{l1n1c}
   \sum_{k=n}^{+\infty} 
\left(\frac{1}{2^k}\right)^{2^*/2-q} 
\int_{\{2^n\leq \theta\leq 2^{k+1}\}}
|\theta|^{N/(N-2)} \dx
\leq C(q) \int_\Omega |\theta|^q \dx,
\end{gather}
where $C(q)$ is a constant only depending on $q$.
\par
After a few computations, from inequality (\ref{l1n1a}) together with
(\ref{l1n1b}) and (\ref{l1n1c}) 
it follows that
\begin{gather} \label{l1n2}
   \int_\Omega |\theta|^q\dx \leq 
M_1(n,q,\Omega,f,C_1) +(\omega(n))^{2^*/2}
M_2\big(\|\theta\|_{L^1(\Omega)},q,C_1,\Omega\big) \int_\Omega |\theta|^q \dx,
\end{gather}
where $M_1$ is a constant only depending on $n$, $q$, $\Omega$, $f$
and $C_1$, and $M_2$ is a constant only depending on
$\|\theta\|_{L^1(\Omega)}$, $q$, $C_1$ and $\Omega$.
\par
Therefore
since  $\omega(n)$ tends to zero as $n$ tends to infinity,
 we can choose $n$ such that $(\omega(n))^{2^*/2}
M_2\big(\|\theta\|_{L^1(\Omega)},q,C_1,\Omega\big)<1/2$
 and then
(\ref{l1n2}) yields
\begin{gather*}
  \int_\Omega |\theta|^q \dx \leq 2M_1(n,q,\Omega,f,C_1),
\end{gather*}
that is (\ref{len2}).

\par\medskip
\noindent{\slshape Case $N=2$.} 
Let $a$ and $M$ be two non negative real numbers such that
\begin{equation}
  \label{l1n4}
  \forall r\in\R^+\quad |f(r)|\leq a+ M|r|.
\end{equation}
\par
Using similar techniques as in the previous case, we obtain
\begin{align*}
  \int_\Omega \Bigg|\frac{D\theta}{1+|\theta|}\Bigg|^2\dx 
& \leq C \left(1+
\sum_{k=0}^{+\infty} \frac{1}{2^{k+1}}\int_{\{0<\theta\leq 2^{k+1}\}}
|\theta|^2\dx \right) \\
& \leq C\left(1+
\int_\Omega |\theta|\dx\right),
\end{align*}
where $C$ is a constant only depending on $\Omega$, $a$, $M$ and $C_1$.
\par
It follows that
\begin{equation}
  \label{l1n5}
  \big\| \ln(1+|\theta|)\big\|_{H^1_0(\Omega)} \leq C,
\end{equation}
where $C$ is a constant only depending on $\Omega$, $a$, $M$, $C_1$
and $\|\theta\|_{L^1(\Omega)}$.
\par
Making use of Theorem 7.15 from \cite{GiTr}, let $C_2$ and $C_3$ be
two positive real numbers only depending on $N$, such that
\begin{equation*}
  \int_\Omega
  \exp\left[
\left(\frac{\ln(1+|\theta|)}{C_2\|D(\ln(1+|\theta|))\|_{L^2(\Omega)}} 
\right)^2\right] \dx \leq C_3\big|\Omega\big|
\end{equation*}
Therefore from (\ref{l1n5}) we have
\begin{equation*}
  \int_\Omega \exp\left[
\left( \frac{\ln(1+|\theta|)}{C_2C} \right)^2
\right] \dx \leq C_3\big|\Omega\big|,
\end{equation*}
where $C$, $C_2$ and $C_3$ does not depend on $\theta$.
\par
If follows that for any $1\leq q <+\infty$ there exists $C'>0$ only
depending on $f$, $\Omega$, $C_1$, $C_2$, $C_3$, $\|\theta\|_{L^1(\Omega)}$ and $q$
such that
\begin{equation*}
  \|\theta\|_{L^q(\Omega)}\leq C'.
\end{equation*}
\par

The proof of Lemma \ref{lem2} is
complete.
\end{proof}

\medskip

We are now in a position to prove Theorem \ref{th2}.
\begin{proof}[Proof of Theorem \ref{th2}]
  The proof is divided into 3 steps. In Step 1 we consider
a solution $(\uep,\thep)$ of the approximate problem
(\ref{eq1})--(\ref{eq3}) with $\fep=f\circ T_{1/\varepsilon}$\ 
($\varepsilon>0$) in place of $f$ and we derive a few preliminary
estimates. In Step 2, using the coupling between the unknowns $\uep$
and $\thep$, we establish an important equality that first implies an
$L^1(\Omega)$--estimate on $\thep$. In Step 3, we make use of Lemma
\ref{lem2} to obtain an $L^2(\Omega)$--estimate on $\fep(\thep)$ and,
at last, we pass to the limit in the approximate problem.
\par\medskip

{\noindent\slshape Step 1.} For $\varepsilon>0$, we consider the
following system
\begin{align}
  \ & \lambda \uep - \diw\big( \A(x)D\uep-\fep(\thep)\big) =
 g\quad  & & \text{ in $\Omega$,} \label{g9} \\
\ &  \mu \thep -\diw\big(\aop(x,D\thep)\big)= \big(\A(x)D\uep-
\fep(\thep)\big)\cdot
D\uep  & & \text{ in $\Omega$,} \label{g10} \\
& \uep=0\qquad\thep=0 & &\text{ on $\partial\Omega$.} \label{g11}
\end{align}

\par
Analyzing the proof of Theorem \ref{th1} allows us to show that
 there exists at least a
weak-renormalized solution $(\uep,\thep)$ of (\ref{g9})--(\ref{g11}).
Indeed as the continuous function $\fep$ is bounded, 
the mapping $\Gamma$ constructed in the proof of Theorem \ref{th1} is
continuous and compact from $L^1(\Omega)$ into a bounded subset of
$L^1(\Omega)$. Then the Schauder fixed--point Theorem allows us to conclude.
\par
For $\varepsilon>0$, let $(\uep,\thep)$ be a weak-renormalized solution of
(\ref{g9})--(\ref{g11}). It follows from (\ref{psm4}) that
\begin{gather}
  \int_\Omega \big|D\uep\big|^2\dx < C \Big( 1+
  \big\|\fep(\thep)\big\|_{(L^2(\Omega))^N}^2 \Big), \label{g12} \\
\big\|(\A(x)D\uep-\fep(\thep))\cdot D\uep\big\|_{L^1(\Omega)} \leq 
C \Big( 1+
  \big\|\fep(\thep)\big\|_{(L^2(\Omega))^N}^2 \Big), \label{g13} 
\end{gather}
and, recalling that $\thep$ is a renormalized solution of
(\ref{g10})--(\ref{g11}) and using Proposition \ref{prop1}, 
\begin{equation}
  \label{g14}
  \int_\Omega \big|DT_K(\thep)\big|^2 \dx \leq C K \Big( 1 + 
\big\|\fep(\thep)\big\|_{(L^2(\Omega))^N}^2 \Big)
\quad\quad \forall K>0,
\end{equation}
where $C$ is a  constant independent of $\varepsilon$ and $K$.

\begin{rmk}\label{ogrmk55}
  From Corollary \ref{coro1} we obtain that for
  any $1\leq p <N/(N-1)$
  \begin{equation}
    \label{g14b}
    \big\|\thep\big\|_{W^{1,p}_0(\Omega)} \leq C \Big(
1+ \big\|\fep(\thep)\big\|_{(L^2(\Omega))^N}^2\Big).
  \end{equation}
In the case where $\lim_{r\rightarrow +\infty} |f(r)|^2/r=0$, 
deriving an $L^q$--estimate
for any $1\leq q <N/(N-2)$ (and $1\leq q <+\infty$ if $N=2$) is an
easy task. But under hypothesis (\ref{g1}), a new estimate on $\thep$
is necessary to obtain a upper bound on $\|\thep\|_{L^1(\Omega)}$ and 
more generally on $\|\thep\|_{L^q(\Omega)}$.
\end{rmk}

\par\medskip
{\noindent \slshape Step 2.} 
For $K>0$, since $T_K(\thep)\in L^\infty(\Omega)\cap H^1_0(\Omega)$
and $DT_K(\thep)=0$ almost everywhere  on $\{x\,:\,|\thep(x)|\geq K\}$ 
 Proposition \ref{prop1} (with
$w=T_K(\thep)/K$) gives that
\begin{gather*}
  {\mu}\int_\Omega \frac{T_K(\thep)}{K}\,\thep \dx +\frac{1}{K}
\int_\Omega \aop(x,D\thep)\cdot
  DT_K(\thep)\dx  = \int_\Omega  \frac{T_K(\thep)}{K}\,
(\A(x)D\uep-\fep(\thep))\cdot D\uep \dx.
\end{gather*}
Plugging the test function $\uep$ in (\ref{g9}) and summing the result
to the previous equality  yield
\begin{multline*}
  \lambda \int_\Omega(\uep)^2\dx+\mu\int_\Omega \frac{T_K(\thep)}{K}\,
  \thep \dx
  +\frac{1}{K}\int_\Omega \aop(x,D\thep)\cdot DT_K(\thep)\dx \\
+  \int_\Omega  \frac{K-T_K(\thep)}{K}\,
\A(x)D\uep\cdot D\uep\dx = 
  \int_\Omega \frac{K-T_K(\thep)}{K}\,
\fep(\thep)\cdot D\uep\dx 
+  \int_\Omega g\uep\dx.
\end{multline*}

Since $ \lambda>0$ and $K-T_K(\thep)$ is a 
non negative function, the coercivity
 of $\aop$ and $\A$ together with Young's inequality lead to
\begin{multline}
  \label{g15b}
   \int_\Omega(\uep)^2\dx+\int_\Omega \frac{T_K(\thep)}{K}\,\thep\dx 
+   \frac{1}{K}\int_\Omega \big| DT_K(\thep)\big|^2\dx \\
+ \int_\Omega\frac{K-T_K(\thep)}{K}\,|D\uep|^2\dx \leq 
C\left(\int_\Omega \frac{K-T_K(\thep)}{K}\,
\big|\fep(\thep)\big|^2\dx+\int_\Omega g^2\dx\right),
\end{multline}
where $C$ is a constant independent of $\varepsilon$ and $K$.
\par
As $\forall \varepsilon>0$
\begin{gather*}
  \thep \frac{T_K(\thep)}{K} \xrightarrow{K\rightarrow 0} |\thep|
  \quad\text{almost everywhere  in $\Omega$ and}\\ 
\frac{K-T_K(\thep)}{K} \xrightarrow{K\rightarrow 0} 2\indi_{\{\thep<0\}} 
+ \indi_{\{\thep=0\}} \quad\text{almost everywhere  in $\Omega$},
\end{gather*}
passing to the limit as $K$ tends to
zero in inequality (\ref{g15b}) gives that, $\forall \varepsilon>0$,
\begin{equation}
  \label{g16}
   \int_\Omega(\uep)^2\dx+\int_\Omega|\thep|\dx
   \leq C\left(\int_\Omega \big|\fep(\thep)\big|^2
 \indi_{\{\thep\leq 0\}} \dx+\int_\Omega g^2\dx\right)
\end{equation}
where $C$ is constant independent of $\varepsilon$.
\par
Due to  assumption (\ref{g2}) on the behavior of $f$ on $\R^-$,
$\forall \eta>0$, $\exists C_\eta>0$ such that, $\forall r\in\R^-$, 
$|f(r)|^2\leq \eta|r|+C_\eta$. If we choose $\eta$ sufficiently small,
then inequality (\ref{g16}) implies that there exists $C_1>0$ such that,
 $\forall\varepsilon>0$,
 \begin{equation}
   \label{g17}
  \int_\Omega(\uep)^2\dx+\int_\Omega|\thep|\dx \leq C_1.  
 \end{equation}

 \begin{rmk}
   Inequality (\ref{g16}) shows that $\thep$ and $\uep$ are controlled
  respectively  in $L^1(\Omega)$ and in $L^2(\Omega)$ if $\fep(\thep)$ is
  controlled in $L^2(\Omega)$ only on the subset of $\Omega$ where
  $\thep\leq0$. Inequality (\ref{g16}) is also used to prove
  uniqueness result (see Theorem \ref{th4}).
 \end{rmk}
\par\medskip
{\noindent \slshape Step 3.} 
As all the  terms in the left hand side of 
(\ref{g15b}) are non negative, one has,  
$\forall \varepsilon>0$ and $\forall K>0$, 
\begin{gather*}
\begin{split}
   \frac{1}{K}\int_\Omega \big| DT_K(\thep)\big|^2\dx 
& \leq
C\left(\int_\Omega \frac{K-T_K(\thep)}{K}\,
\big|\fep(\thep)\big|^2\dx+\int_\Omega g^2\dx\right) \\
& \leq C \left(\int_{\{\thep< 0\}} 2\big|\fep(\thep)\big|^2\dx
+ \int_{\{0\leq\thep<K\}} \big|\fep(\thep)\big|^2\dx +
\int_\Omega g^2\dx\right),
\end{split}
\end{gather*}
since $K-T_K(\thep)=0$ almost everywhere  on $\{x\,:\,\thep(x)\geq K \}$.
\par
Therefore  growth assumption (\ref{g1}) on $f$  and estimate (\ref{g17})
 imply  that there exists $C_2>0$ such
that  $\forall \varepsilon>0$,  $\forall K>0$
\begin{equation*}
  \frac{1}{K}\int_\Omega \big| DT_K(\thep)\big|^2\dx \leq
C_2\left(\int_{\{0\leq \thep\leq K\}} \big|\fep(\thep)\big|^2
  \dx+ 1 \right).
\end{equation*}
Let us denote $f^*$ the real-valued function defined by  $f^*(r)=\sup_{0\leq r'\leq
  r}|f(r')|$, $\forall r\in\R^+$. The function $f^*$ satisfies (\ref{g1}) and
  $\forall\varepsilon>0$, $\forall K>0$
\begin{equation*}
  \frac{1}{K}\int_\Omega \big| DT_K(\thep)\big|^2\dx \leq
C_2\left(\int_{\{0\leq \thep < K\}} \big|f^*(\thep)\big|^2
  \dx+ 1 \right).
\end{equation*}
\par
Since $C_1$ and $C_2$ are independent of $\varepsilon$ and $K$,
from (\ref{g16}) and  the above
inequality we can apply 
 Lemma \ref{lem2} to $\thep$, $\forall \varepsilon>0$. It follows that
  the sequence $\thep$ is bounded in $L^q(\Omega)$
for any $1\leq q <N/(N-2)$ (and $1\leq q <+\infty$ if
$N=2$).
In particular, growth  assumptions (\ref{g1}) and (\ref{g2}) on $f$
 imply that
\begin{equation}
  \label{g18}
\fep(\thep) \text{ is bounded in }\big(L^2(\Omega)\big)^N,
\end{equation}
and from (\ref{g14b})  we obtain that for any $1\leq p <N/(N-1)$ 
\begin{equation*}
  \thep \text{ is bounded in } W^{1,p}_0(\Omega).
\end{equation*}
By the Rellich Kondrachov Theorem, let $\theta$ be a measurable
  function defined from $\Omega$ into $\R$ such that, up to a
  subsequence, 
  $\forall 1\leq q< N/(N-2)$ (and $1<q<+\infty$ if $N=2$)
\begin{equation}
\quad \thep\longrightarrow \theta \text{ in
  $L^q(\Omega)$ and almost everywhere  in $\Omega$, as $\varepsilon$ tends to zero.} 
  \label{g19}
\end{equation}
\par
Since $(N+2)/N<N/(N-2)$, the continuity of $f$, 
growth assumptions (\ref{g1}) and (\ref{g2}) and (\ref{g19}) allow us to
deduce,  by a
standard equiintegrability argument, that
\begin{equation}
  \label{g22}
  \fep(\thep) \longrightarrow f(\theta)\quad\text{in
  $\big(L^2(\Omega)\big)^N$ as $\varepsilon$ tends to zero.}
\end{equation}
\par
Next, using the linear character of equation (\ref{g9})  with respect to
$\uep$ together with (\ref{g22}) it is easy to show that
\begin{equation*}
  \uep \longrightarrow u \quad\text{in $H^1_0(\Omega)$ as
  $\varepsilon$ tends to zero,}
\end{equation*}
and then, $(u,\theta)$ satisfies equation (\ref{eq1}) in ${\mathcal
  D'(\Omega)}$ with $u\in
H^1_0(\Omega)$ and $f(\theta)\in\big(L^2(\Omega))^N$.
\par
 It follows that
\begin{gather*}
 \big(\A(x)D\uep-\fep(\thep)\big)\cdot D\uep  
\longrightarrow \big(\A(x)Du-f(\theta)\big)\cdot Du 
\quad\text{in $L^1(\Omega)$ as $\varepsilon$ tends to zero.}
\end{gather*}
\par
As far as equation (\ref{g10}) is concerned, the 
 continuity of renormalized solution with respect to the
data implies that $\theta$ is a renormalized solution of
 (\ref{eq2}).\par
The proof of Theorem \ref{th2} is complete.
\end{proof}

\medskip
We now prove Theorem \ref{th3}.

\begin{proof}[Proof of Theorem \ref{th3}]
Let $\ftid$ be the function defined by
\begin{equation*}
  \ftid(r)=\left\{
    \begin{aligned}
      0 & \text{ if }r\leq r_0, \\
f(r) & \text{ if } r>r_0.
    \end{aligned}
\right.
\end{equation*}

The function 
 $\ftid$ is continuous and satisfies assumptions  (\ref{g1}) and (\ref{g2}). 
Making use of Theorem \ref{th2}, let $(u,\theta)$ be a
 weak-renormalized solution of system (\ref{eq1})--(\ref{eq3}) with
 $\ftid$ in place of $f$.
\par
Our aim now is to prove that
\begin{equation*}
  \theta\geq r_0\quad\text{almost everywhere  in $\Omega$.}
\end{equation*}
\par
For $K>0$, let $H$ be the function defined by $H(r)=-T_K^-(r-r_0)$,
$\forall r\in\R$.
We  have $H\in W^{1,\infty}(\R)$,  $H'(r)=\indi_{\{-K+r_0<r<r_0\}}$,
 so $H'$ has a compact support.
Since $r_0\leq 0$, it follows that $H(\theta)\in
L^{\infty}(\Omega)\cap H^1_0(\Omega)$ and recalling that $\theta$ is a 
renormalized solution of (\ref{eq2}),
Proposition \ref{prop1} with $w=H(\theta)$ leads to
\begin{gather*}
  \mu\int_\Omega \theta H(\theta)\dx + \int_{\{-K+r_0<\theta<r_0\}}
 \aop(x,D\theta)\cdot D\theta \dx
=-\int_\Omega \big(\A(x)Du-\ftid(\theta)\big)\cdot Du H(\theta)\dx.
\end{gather*}
The definitions of $H$ and $\ftid$ imply that $\ftid(r) H(r)=0$,
$\forall r\in\R$, and because $H(r)\leq 0 $
 the
coercivity of $\aop$ and $\A$ gives
\begin{equation*}
  \int_\Omega |\theta|\,T_K^-(\theta-r_0)\dx \leq 0.
\end{equation*}
\par
It follows that
\begin{equation*}
  \theta\geq r_0\quad\text{almost everywhere  in $\Omega$,}
\end{equation*}
and according to the definition of $\ftid$,
\begin{equation*}
  \ftid(\theta)=f(\theta)\quad\text{almost everywhere  in $\Omega$.}
\end{equation*}
\par
Hence $(u,\theta)$ is a weak-renormalized solution of (\ref{eq1})--(\ref{eq3}).
\end{proof}

\medskip
\section{Uniqueness results}

In this section we  assume that
\begin{equation}
  \label{un1}
  f(0)=0,
\end{equation}
and we give the following uniqueness result of a
small  solution  $(u,\theta)$ of (\ref{eq1})--(\ref{eq3}) 
such that $\theta\geq 0$ almost everywhere  in $\Omega$ under
additional assumptions on $f$, $\aop$ and $N$.

\par\smallskip

\begin{theorem}\label{th4}
  Assume that assumptions (\ref{as1})--(\ref{as5}), 
  (\ref{g1}) and (\ref{un1}) hold true. Moreover assume that 
  \begin{gather}
   N=2 \quad \text{or}\quad N=3, \label{th4h1} \\
    \big(\aop(x,\xi)-\aop(x,\xi')\big)\cdot(\xi-\xi') \geq
    \delta|\xi-\xi'|^2 \quad\text{almost everywhere  in $\Omega$, $\forall
    \xi,\xi'\in\R^N$,} \label{th4h2}\\
\exists L>0\text{ such that } \forall r,r'\in\R^+\quad |f(r)-f(r')|\leq
    L|r-r'|.  \label{th4h3}
  \end{gather}
There exists $\eta>0$ such that if $\|g\|_{L^2(\Omega)}<\eta$, then
the weak-renormalized solution $(u,\theta)$ 
of (\ref{eq1})--(\ref{eq3}), such that
$\theta\geq0$ almost everywhere  in $\Omega$, is unique.
\end{theorem}
\par\smallskip

\begin{proof}[Proof of Theorem \ref{th4}]
From Theorem \ref{th3} let $(u_1,\theta_1)$ and $(u_2,\theta_2)$ be
two weak-renormalized solutions of (\ref{eq1})--(\ref{eq3}) such that
$\theta_1\geq 0$ and $\theta_2\geq 0$ almost everywhere  in $\Omega$.
\par
The aim is to prove that
\begin{equation}
  \label{t4n1}
  \big\|\theta_1-\theta_2\big\|_{L^2(\Omega)} \leq
  \omega\big(\|g\|_{L^2(\Omega)}\big)
  \big\|\theta_1-\theta_2\big\|_{L^2(\Omega)}, 
\end{equation}
where $\omega$ is independent of $\theta_1$ and $\theta_2$ and is such
that $\omega(r)$ tends to zero as $r$ tends to zero. 
\par
We denote by $F_i$ the term $\big(\A(x)Du_i-f(\theta_i))\cdot Du_i$,
for $i=1,2$. Proposition \ref{prop2} and (\ref{th4h2}) give 
\begin{equation*}
  \forall K>0,\quad 
\delta \int_{\{|\theta_1-\theta_2|<K\}} \big| D\theta_1-D\theta_2\big|^2 \dx
\leq K \int_\Omega \big|F_1-F_2\big| \dx.
\end{equation*}
From a result of \cite{DMOP99}, it follows that
$T_K(\theta_1-\theta_2)\in H^1_0(\Omega)$ for any $K>0$. As 
$N=2$ or $N=3$ there exists $1<p<N/(N-1)$ such that $p^*=2$ and so
Lemma \ref{lem1} and the above inequality imply that
\begin{equation}
  \label{t4n2}
  \big\|\theta_1-\theta_2\big\|_{L^2(\Omega)} \leq C
  \big\|F_1-F_2\big\|_{L^1(\Omega)}, 
\end{equation}
where $C$ is a generic constant independent of $i$ and $g$.
\par
A  calculus leads to
\begin{align*}
    \big\|F_1-F_2\big\|_{L^1(\Omega)} \leq & 
    \|A\|_{(L^\infty(\Omega))^{N\times N}}
\big\|Du_1-Du_2\big\|_{(L^2(\Omega))^N} \times 
\big\|Du_1+Du_2\big\|_{(L^2(\Omega))^N} \\
& {} + \big\|f(\theta_1)-f(\theta_2)\big\|_{(L^2(\Omega))^N}\times
    \big\|Du_1\big\|_{(L^2(\Omega))^N}
\\ & {} +
 \big\|f(\theta_1)\big\|_{(L^2(\Omega))^N}\times 
\big\|Du_1-Du_2\big\|_{(L^2(\Omega))^N}.
\end{align*}
The linear character of equation (\ref{eq1}) gives
\begin{gather*}
  \|Du_1-Du_2\|_{(L^2(\Omega))^N} \leq C
  \big\|f(\theta_1)-f(\theta_2)\big\|_{(L^2(\Omega))^N}, \\
\|Du_i\|_{(L^2(\Omega))^N} \leq C \big(
  \|f(\theta_i)\|_{(L^2(\Omega))^N} +
  \|g\|_{L^2(\Omega)}\big),\quad\text{for $i=1,2$}.
\end{gather*}
Using (\ref{th4h3}) and the above inequalities we obtain
\begin{equation}
  \label{t4n3}
  \big\|F_1-F_2\big\|_{L^1(\Omega)} \leq 
C \big( \|f(\theta_1)\|_{(L^2(\Omega))^N} +
 \|f(\theta_2)\|_{(L^2(\Omega))^N} + \|g\|_{L^2(\Omega)}\big)\,
\big\|\theta_1-\theta_2\big\|_{L^2(\Omega)},
\end{equation}
and therefore (\ref{t4n1}) reduces to prove that
\begin{equation}
  \label{t4n4}
  \big\|f(\theta_i)\big\|_{(L^2(\Omega))^N} \leq
  \omega\big(\|g\|_{L^2(\Omega)}\big),\quad \text{for $i=1,2$,}
\end{equation}
 where $\omega$ is independent of $i$ and is such
that $\omega(r)$ tends to zero as $r$ tends to zero. 
\par
Since $\theta_i\geq 0$ almost everywhere  in $\Omega$, (\ref{g16}) 
implies that
\begin{equation}
  \label{t4n5}
  \big\|\theta_i\big\|_{L^1(\Omega)} \leq C\|g\|_{L^2(\Omega)}\quad
  \text{for $i=1,2$,}
\end{equation}
and from (\ref{g15b}) we obtain, for $i=1,2$,
\begin{equation}
  \label{t4n6}
  \forall K>0\quad\frac{1}{K} \int_\Omega \big|DT_K(\theta_i)\big|^2
  \dx \leq C\left( \int_{\{0\leq \theta_i \leq K\}}
  \big|f(\theta_i)\big|^2 \dx + \|g\|_{L^2(\Omega)}^2\right).
\end{equation}
So if $\|g\|_{L^2(\Omega)}\leq 1$,  (\ref{t4n5}) and 
(\ref{t4n6}) together with Lemma \ref{lem2} allow us to deduce that, $\forall 1\leq q <N/(N-2)$ (and
$q<+\infty$ if $N=2$),
\begin{equation}
  \label{t4n7}
  \big\|\theta_i\big\|_{L^q(\Omega)} \leq C(q),
\end{equation}
where $C(q)$ is a constant independent of $i$ and $g$.
\par
Using (\ref{g1}), (\ref{un1}) and (\ref{th4h3}), let $M>0$ such that 
\begin{equation*}
  \forall r\in\R^+\quad |f(r)|\leq M r^{(N+2)/(2N)}.
\end{equation*}
By interpolation between $L^1(\Omega)$ and $L^{(N+1)/(N-1)}(\Omega)$
we have, for $i=1,2$,
\begin{equation*}
  \big\|f(\theta_i)\big\|_{(L^2(\Omega))^N} \leq M
  \Big(\big\|\theta_i\big\|_{L^\frac{N+2}{N}(\Omega)}\Big)^\frac{N}{2(N+2)} 
\leq M \Big(\|\theta_i\|_{L^1(\Omega)}\Big)^\frac{1}{2N} 
\Big(\|\theta_i\|_{L^\frac{N+1}{N-1}(\Omega)}\Big)^\frac{N+1}{2N},
\end{equation*}
and using (\ref{t4n5}) and (\ref{t4n7}) 
(indeed $\frac{N+1}{N-1}<\frac{N}{N-2}$), if
$\|g\|_{L^2(\Omega)}\leq 1$ then we have
\begin{equation}
\label{t4n10}
  \big\|f(\theta_i)\big\|_{(L^2(\Omega))^N} \leq C
  \|g\|_{L^2(\Omega)}^\frac{1}{2N},
\end{equation}
where $C$ is independent of $i$ and $g$.
\par
It follows from (\ref{t4n2}), (\ref{t4n3}) and (\ref{t4n10}) that 
(\ref{t4n1}) is proved
for  $\|g\|_{L^2(\Omega)}\leq 1$.
Then there exists $\eta>0$ such that if 
 $\|g\|_{L^2(\Omega)}< \eta$ then
$\omega\big(\|g\|_{L^2(\Omega)}\big)<1$ and (\ref{t4n1}) implies that
$\theta_1=\theta_2$ almost everywhere  in $\Omega$.
\par
The proof of Theorem \ref{th4} is complete.
\end{proof}
\par
\medskip

\bigskip

\noindent{\small\slshape Acknowledgement---}\small The author thanks
D. Blanchard and  F. Murat
for numerous discussions on the problems under investigation in the
present paper.

\par\medskip

\bibliographystyle{plain}
\makeatletter
\renewenvironment{thebibliography}[1]{%
  \@xp\section\@xp*\@xp{\refname}%
  \normalfont\labelsep .5em\relax
  \renewcommand\theenumiv{\arabic{enumiv}}\let\p@enumiv\@empty
  \list{\@biblabel{\theenumiv}}{\settowidth\labelwidth{\@biblabel{#1}}%
    \leftmargin\labelwidth \advance\leftmargin\labelsep
    \usecounter{enumiv}}%
 \sloppy \clubpenalty\@M \widowpenalty\clubpenalty
  \sfcode`\.=\@m
}
\makeatother
\renewcommand{\baselinestretch}{1.1}\normalsize
\bibliography{artog}

\end{document}